\documentclass[12pt]{article}
\usepackage{amsmath,amsfonts}
\usepackage{e-jc}
\usepackage{amssymb}
\usepackage{appendix}

\usepackage{color}
\newcommand{\refm [1]} {(\ref{#1})}

\newcounter{theorem}

\newtheorem{proposition}{Proposition}
\newtheorem{theorem}{Theorem}
\newtheorem{lemma}{Lemma}

\newcommand{\la}{\label}

\newcommand{\PP}{\mathbb{P}}

\newcommand{\nicef}{{\bf {\cal F}}}
\newcommand{\en}{\end{equation}}

\def\PP{{\mathbb P}}

\def\BR{{\mathbb R}}

\def\BE{{\mathbb E}}

\def\BZ{{\mathbb Z}}
\def\BC{{\mathBbb C}}

\def\proof{\noindent{\bf Proof\ \ }}
\def\qed{\mbox{\rule{0.5em}{0.5em}}}

\def\la{\label}

\def\del1{{\delta_n^{(1)}}}

\newcommand{\ignore}[1]{}

\def\BR{{\mathbb R}}

\def\BE{{\mathbb E}}
\def\BZ{{\mathbb Z}}
\def\BC{{\mathbb C}}
\def\calh{{\cal H}}

\newcommand{\be}{\begin{equation}}
\newcommand{\eu}{\end{equation}}
\newcommand{\ber}{\begin{eqnarray}}
\newcommand{\ena}{\end{eqnarray}}
\newcommand{\nin}{\noindent}
\newcommand{\non}{\nonumber}

\title{Developments in the Khintchine-Meinardus probabilistic method
for asymptotic enumeration}

\author{ Boris L. Granovsky\\
\small Department of Mathematics\\[-0.8ex]
\small Technion-Israel Institute of Technology\\[-0.8ex]
\small Haifa, 32000, Israel\\
\small\tt mar18aa@techunix.technion.ac.il
\and  
Dudley Stark\\
\small School of Mathematical Sciences\\[-0.8ex]
\small Queen Mary, University of London\\[-0.8ex]
\small London E1 4NS, United Kingdom\\
\small\tt D.S.Stark@maths.qmul.ac.uk
}

\date{\dateline{August 7, 2014}{November 4, 2015}}

\begin{document}

\maketitle

\begin{abstract}
A theorem of Meinardus provides asymptotics of the number of weighted
partitions under certain assumptions on associated ordinary and Dirichlet
generating functions.
The ordinary generating functions are closely related to
Euler's generating function
$\prod_{k=1}^\infty S(z^k)$ for partitions, where $S(z)=(1-z)^{-1}$.
By applying a method due to Khintchine,
we extend Meinardus' theorem to find the asymptotics of
the Taylor coefficients of generating functions of the form
$\prod_{k=1}^\infty S(a_kz^k)^{b_k}$ for sequences $a_k$, $b_k$
and general $S(z)$.
We also reformulate the hypotheses of the theorem in terms of
the above generating functions.
This allows novel applications of the method. In particular, we prove rigorously the asymptotics of Gentile statistics and derive the asymptotics of combinatorial objects with distinct components.


\end{abstract}

\nin \section{Introduction}
Meinardus \cite{M} proved a theorem about the asymptotics of weighted
partitions with weights satisfying certain conditions.
His result was extended to the combinatorial objects
called assemblies and selections in \cite{GSE} and to
Dirichlet generating functions for weights, with multiple singularities in
\cite{GS}. In this paper, we extend Meinardus' theorem further to a general
framework, which encompasses a variety of models in physics and combinatorics, including previous results.

Let $f$ be a
generating function of a nonnegative sequence $\{c_n,\ n\ge 0,\ c_0=1\}$:
\be
f(z)=\sum_{n\ge 0}c_n z^n,\la{sac}
\en
with radius of convergence 1. As an example,
consider the number of weighted partitions $c_n$ of size $n$,
 determined by the generating function identity
\be\label{intro}
\sum_{n=0}^\infty c_nz^n=\prod_{k=1}^\infty (1-z^k)^{-b_k}, \quad \vert z\vert<1,
\en
for some sequence of real numbers $b_k\ge 0,\ k\ge 1$.
If $b_k=1$ for all $k\ge 1$, then $c_n$ is the number of integer partitions.
Meinardus \cite{M} proved a theorem giving the asymptotics of $c_n$
under certain assumptions on the sequence $\{b_k\}$.

The generating function in \refm[intro] may
be expressed as $\prod_{k=1}^\infty \big(S(z^k)\big)^{b_k}$, where
$S(z)=(1-z)^{-1}$.
This observation allows the following natural generalization.
Let $f$ in \refm[sac] be of the form:\be\label{frame}
f(z)=\prod_{k=1}^\infty \big(S(a_k z^k)\big)^{b_k},
\en
 with given sequences $0< a_k\le 1$, $ b_k\ge  0,\ k\ge 1$, and a given function $S(z).$

This is a particular case of the  class of general multiplicative models, introduced and studied by Vershik (\cite{V1}). In the setting \refm[frame], in the case of weighted partitions,  a combinatorial meaning can be attributed
to the parameters  $a_k, b_k$. Namely, if $b_k=1,$ then $a_k$ can be viewed  as a properly scaled number of colours for each component of size $k,$ such that
given $l$ components of size $k$, the total number of colourings is
$a_k^l$. Equivalently, to each particular partition of $n,$ say $n=\sum_{l=1}^n j_lk_l,$ is prescribed the  weight $\prod_{l=1}^n a_{k_l}^{j_l}.$  On the other hand, if $a_k=1$, then given $l$ components of size $k$,
%$\binom b_k+l-1 \choose b_k $
the total number of colourings equals  the number of distributions of $b_k$ indistinguishable balls among $l$ cells, so that in this model  $b_k$ has    a meaning of a scaled number of types  prescribed to a component of size $k$.

Yakubovich (\cite{YA}) derived the limit shapes for models \refm[frame] in the case $a_k=1,\ k\ge 1$, under some analytical conditions on  $S$ and $b_k$.
Note that %as has been the case in
past versions \cite{GSE}--\cite{GS} of Meinardus' theorem deal with the asymptotics of $c_n,\ n\to \infty,$ when $a_k=1,\ k\ge 1,$ for three  cases of the function $S,$ corresponding to the three classic models of ideal gas in  statistical mechanics, namely, Maxwell-Boltzmann, Bose-Einstein and Fermi-Dirac. They are mathematically equivalent to aforementioned combinatorial models: assemblies, weighted partitions and selections, respectively. In this context, $n$ is a total energy of a system, which is partitioned into integer energy levels $k_1,\ldots, k_n,$ so that $n= j_1k_1+ \ldots+ j_nk_n,$ where $j_1,\ldots,j_n$ are occupation numbers which are equal to the numbers of particles at corresponding energy levels.

 Our objective in this paper is  to derive the  asymptotics $c_n,\ n\to \infty,$ in the general framework \refm[frame].
 %We will assume that the conditions $(I)- (III)$ below are satisfied.
We suppose that $S(0)=1$ and
 that $S(z)$ can be expanded in a power series with radius of convergence $\ge 1$ and non-negative coefficients $d_j$:
 \be\label{sz}
S(z)=\sum_{j=0}^\infty d_jz^j,
 \en
with $d_0=1$,
  and that
 $\log S(z)$ can be expanded as
\be\label{logSexp}
\log S(z)=\sum_{j=1}^\infty \xi_j z^j,
\en
with  radius of convergence $1$.
We note that the assumption $d_j\ge 0,\ j\ge 0$ is necessary for the implementation of the Khintchine probabilistic approach (see
\refm[rep],\refm[phidef] below).

From \refm[frame] and \refm[logSexp] one can  express the coefficients $\Lambda_k$ of the power series for the function $\log f(z)$, with radius of convergence 1:
\be\label{logflam}
\log f(z) = \sum_{k=1}^\infty \Lambda_k z^k,\quad\Lambda_k=\sum_{j\mid k}b_j a_j^{k/j} \xi_{k/j}.
\en
We define the Dirichlet generating function for the sequence $\Lambda_k:$
\be\label{Dirdef}
D(s)=\sum_{k=1}^\infty \Lambda_k k^{-s},
\en
which by virtue of   \refm[logflam] admits the following representation
\be \label{dpres}
D(s)=\sum_{k=1}^\infty\sum_{j=1}^\infty b_k \xi_j a_k^j
(jk)^{-s},
\en
 as long as $\Re(s)$ is large enough so that the double Dirichlet series in \refm[dpres]
converges absolutely.
If
 $a_k=a,\ 0<a\le 1 $ for all $k\geq 1$, %as has been the case in
%past versions of Meinardus' theorem,
then
$D(s)$ can be factored as
\be\label{factor}
D(s)=D_b(s)D_{\xi,a}(s),
\en
where
$$
D_{\xi,a}(s)=\sum_{j=1}^\infty a^j\xi_j j^{-s}
$$
and
\be\label{Dbdef}
D_b(s)=\sum_{k=1}^\infty b_k k^{-s}.
\en

The greater generality of \refm[frame] than in previous versions of Meinardus'
theorem will allow novel applications. The proof of
Theorem~\ref{general}, stated below, is a substantial modification
of the method used in
\cite{GSE,GS2,GS}.

We suppose that $\Lambda_k$ and $D(s)$ satisfy conditions $(I)- (III)$,
which are modifications of the three original Meinardus' conditions in \cite{M}.
\bigskip

{\bf Condition $(I)$.}
The Dirichlet generating function $D(s), \ s=\sigma+it$ is analytic in the half-plane $\sigma>\rho_r>0$ and
it  has $r\geq 1$ simple poles at positions $0<\rho_1<\rho_2<\ldots<\rho_r,$ with positive residues
$A_1,A_2,\ldots,A_r$ respectively.
It may also happen that $D(s)$ has a simple pole at $0$ with residue $A_0$.
(If $D(s)$ is analytic at 0, we take
$A_0=0$). Moreover,
 there is a constant $0<C_0\le 1,$ such
that the function $D(s)$, $s=\sigma+it$, has a meromorphic continuation
 to the half-plane
$$
\calh=\{s:\sigma\geq-C_0\}
$$
on which it is analytic except for
the above $r$ or $r+1$ simple poles.

{\bf Condition $(II)$.}
 There is a constant $C_1>0$ such that
$$
D(s)=O\left(|t|^{C_1}\right),\quad t\to\infty
$$
 uniformly for $s=\sigma+it\in\calh$.

{\bf Condition $(III).$}

  The following property of the parameters $a_k$, $b_k$ holds:
\be\label{assum1}
b_ka_k^{l_0} \geq C_2 k^{\rho_r-1},\ k\ge 1,\ C_2>0,
\en
where
\be\label{l0def}
l_0:=\min\{j>0:d_j>0\}.
\en

Moreover, if $l_0>1$ then
%for  $\delta_n$ as defined below in \refm[delasy
   for some fixed $\epsilon>0$ and for small enough $\delta$,
\be\la{III}
2\sum_{k=1}^{\infty}
 {\Lambda_k{e^{-k\delta}\sin^2(\pi k\alpha)}}
\ge \left(1+\frac{\rho_r}{2}+\epsilon\right)
|\log \delta|,\ \
 (2 l_0)^{-1}\leq
|\alpha|\leq 1/2, \ \  l_0>1,
\en
where $\Lambda_k$ is as defined in \refm[logflam].

In order to state our main result, we need some more notations,
which were also used in \cite{GS}.
Define the finite set
\be\la{sd1}
\tilde{\Upsilon}_r=
\left\{\sum_{k=0}^{r-1} \tilde{d}_k(\rho_r-\rho_k):\   \tilde{d}_k\in\BZ_+,\ \sum_{k=0}^{r-1}\tilde{d}_k\geq 2\right\}\cap \big(0,\rho_r+1\big],\quad r\ge 1,
\en
where we have set $\rho_0=0$ and let  $\BZ_+$ denote the set of nonnegative integers. Let $0<\alpha_1<\alpha_2<\cdots<\alpha_{|\tilde{\Upsilon}_r|}\le \rho_r+1$ be all ordered numbers forming the
set $\tilde{\Upsilon}_r$. Clearly, $\alpha_1=2(\rho_r-\rho_{r-1}),$ if the set $\tilde{\Upsilon}_r$ is not empty.
We also
 define  the finite set
\be\la{sd2} \Upsilon_r=\tilde{\Upsilon}_r\cup\{\rho_r-\rho_k:\ k=0,1,\ldots,r-1\},
\en
observing that  some of the differences $\rho_r-\rho_k,\ k=0,\dots,r-1$ may fall into the set $\tilde{\Upsilon}_r$. We let $0<\lambda_1<\lambda_2<\cdots<\lambda_{|\Upsilon_r|}$
be all ordered numbers forming the
set $\Upsilon_r$.

%In addition to its appearance in (\ref{sd2}), the set  $\tilde{\Upsilon}_r$ plays an auxiliary role in the proof of Proposition~\ref{deltaexp} below.

%{\bf Remark} Let $a_k=a$, $0<a<1$, $k\ge 1$ and let the condition $(I)$ hold.  By the root test applied to \refm[logSexp], $\limsup_{j\to \infty} (\vert \xi_j\vert)^{1/j}=1.$  Viewing \refm[Dbdef] as a power series in $a$, we again imply the root test to conclude that the series $D_{\xi}(s)$ converges absolutely and is entire for all $s\in \nicec$.  Hence, by \refm[factor] and \refm[Dbdef], all $r$ poles in condition $(I)$, belong to the function $D_b(s)$, and we may derive from Wiener-Ikehara theorem (see Theorem 2.2, p.122 in \cite{K}) the  upper bound $b_k=o(k^{\rho_r}),\ k\to \infty$, on the rate of growth of $b_k$ (for more details see \cite{GS2}).  With no restrictions on $a_k$, Condition $(I)$ does not allow any conclusion about the poles of $D_b$. This can be easily seen when $a_k=1,\ k\ge 1$.

\begin{theorem}\label{general}
Suppose conditions $(I) - (III)$ are satisfied.

Suppose that $c_n$ has ordinary generating function of the form
\refm[frame], where $0< a_k\le 1$ and $b_k\ge 0, \ k\ge 1$,
that \refm[assum1] in Condition $(III)$ is satisfied
for a constant $C_2>0$,
and that
\be\label{assum2}
\frac{d^2}{d\delta^2}\log S\left(e^{-\delta}\right)>0,
\quad\delta>0.
\en
We then have, as $n\to\infty$,
\be\label{cnasymp}
c_n\sim Hn^{-\frac{2+\rho_r-2A_0}{2(\rho_r+1)}}\exp \Big(\sum_{l=0}^{r}P_l\, n^{\frac{\rho_l}{\rho_r+1}}+
\sum_{l=0}^r\hat{ h}_l\sum_{s:\lambda_s\le \rho_l}K_{s,l}\, n^{\frac{\rho_l-\lambda_s}{\rho_r+1}} \Big),
\en
where $H$, $P_l$, $\hat{h}_l$ and $K_{s,l}$ are constants.
In particular, if  $r=1$, then
$K_{s,l}=0$ for all $s$ and $l$,
$$P_1=\left(1+\frac{1}{\rho_1}\right)
\big(A_1\Gamma(\rho_1+1)\big)^{1/(\rho_1+1)}
$$
and
$$
H=e^{\Theta-\gamma A_0}
\left(2\pi(1+\rho_1)\right)^{-1/2}(A_1\Gamma(\rho_1+1))^{\frac{1-2A_0}{2(\rho_1+1)}},
$$
where
\be
\Theta:=\lim_{s\to 0}(D(s)-A_0s^{-1})
\la{THETA}\eu
and $\gamma$ is Euler's constant.
\end{theorem}
{\bf Remark} The sums in \refm[cnasymp] could be taken from $l=1$ to $r$
and $P_0$ and $\hat{h}_0K_{0,0}$ absorbed into $H$,
but we prefer not to do that as the constants
in \refm[cnasymp] arise naturally from the proof of Theorem~\ref{general}.
The constants $P_l$, $K_{s,l}$ in \refm[cnasymp] are calculated by the
recursive method of \cite{GS}. We do not repeat the description of the method
here.

Theorem~1 generalizes the seminal results by Khintchine \cite{Kh} and Meinardus \cite{M}, as well as their extensions  in \cite{GSE,GS}, and implies
the results therein, including expansive weighted partitions,
 for which
$S(z)=(1-z)^{-1}$, $a_k=1,\ k\ge 1$ and $b_k=k^{r-1},\ k\ge 1$ for some $r>0$.

{\bf Example} This example shows that
\refm[assum1] is  not implied by the other hypotheses of
Theorem~\ref{general}.
Let $a_k=1$ for all $k$, let $b_k=k^{\rho_1-1}$
 and let $\xi_k=k^{\rho_2-1}$, where $0<\rho_1<\rho_2$.
Then, $D_{\xi,1}(s)=\zeta(s+1-\rho_2)$, $D_b(s)=\zeta(s+1-\rho_1)$,
where $\zeta$ is the Riemann zeta function, and $D(s)=D_b(s)D_{\xi,1}(s)$
has simple poles at $\rho_1$ and $\rho_2$. Moreover,
$S(z)=\exp\left(\sum_{k=1}^\infty k^{\rho_2-1}z^k\right)$ has radius of convergence
$1$ and it is easy to check that \refm[assum2] is satisfied.
However,\refm[assum1] is violated.
%Theorem~\ref{general} is proven in Section~2.  In the remaining two sections, we focus on two novel applications implied by  Theorem~\ref{general}. In  Sections~3 and 4 we apply our results to the asymptotic enumeration of Gentile statistics and expansive selections with $a_k=k^{-q}.$ The latter generalizes previous results for polynomials over a finite field.
%\section{Proof of Theorem ~\ref{general}}

As in \cite{GSE}-\cite{GS}, the proof of Theorem~\ref{general} is based on the  Khintchine type representation \cite{Kh}
\be
c_n=e^{n\delta} f_n(e^{-\delta})\PP\left(U_n=n\right),\quad n\ge
1,
 \la{rep}\eu
where $\delta>0$ is a free parameter,
\be
f_n(z)=\prod_{k=1}^n S(a_k z^k)^{b_k}
\la{fndef}\eu
is the $n$-truncation of $f$ in \refm[frame], and
$U_n,\ n\ge 1$ are integer-valued
random variables with characteristic functions
defined by
\be\label{phidef}
\phi_n(\alpha)=\BE\left(e^{2\pi i\alpha U_n}\right)=
\prod_{k=1}^n
\left(\frac{S\left(a_k e^{2\pi i k\alpha-k\delta}\right)}{S\left(a_ke^{-k\delta}\right)}\right)^{b_k},\quad n\ge 1,\quad \alpha\in \BR.
\en
Khintchine established \refm[rep] for the three basic models of statistical mechanics. For general multiplicative measures \refm[rep] was stated in
equation (4) of \cite{GSE}.
It remains to analyze each of the three factors of the right hand side of \refm[rep].

The proof of Theorem~\ref{general} is similar in form to proofs
in \cite{GSE,GS2,GS}, however there are notable differences.
The form of the generating function $f(z)$ in
\refm[frame] is much more general than in the cases of the aforementioned classic combinatorial structures(=models of statistical mechanics). Also recall that in the  previous papers,
the parameters $a_k$ in \refm[frame] were always taken to
be equal 1.
However, we still have
a nice representation of $D(s)$ given in \refm[factor], which allows us to
proceed with the Meinardus-Khintchine method.
The basic method of proof is an analysis of the three factors of
\refm[rep] when $\delta=\delta_n$ is chosen to be the solution of
the equation $\BE U_n=n$, $n\geq 1$. The convexity assumption \refm[assum2]
is made in order to guarantee that there exists a unique solution $\delta_n$ whose
asymptotics we may obtain.
Regarding conditions $(I)$, $(II)$, and $(III)$,
conditions $(I)$ and $(II)$ are similar to the corresponding conditions
in \cite{GS}, in which Dirichlet generating functions with multiple poles
were considered.
An asymptotic equation for $\delta_n$ is obtained from condition $(I)$
and $(II)$ by using Mellin transforms
and then changing the contour of integration.
Assumption \refm[III] of Condition $(III)$ is stronger than the
original condition of Meinardus, however here it is used to prove a local limit
theorem, the proof of which requires, in addition, assumption \refm[assum1],
which involves both $a_k$ and $b_k$.

The proof of Theorem~\ref{general} is
contained in Appendix~\ref{proofs}.
\bigskip

{\bf Two historical remarks} 
\bigskip

1. Khintchine's
probabilistic approach for asymptotic enumeration resulted in
the representation \refm[rep].
Explaining his idea to replace the saddle point method with
a local limit theorem, Khintchine wrote in \cite{Kh}:".......the main novelty of this approach consists of replacing the
complicated analytical apparatus (the method of Darwin-Fowler) by
...the well developed limit theorems of the theory of
probability...that can form the analytical basis for all the
computational formulas of statistical mechanics." In \cite{GSE}, the practical advantages of this idea regarding the asymptotic problems considered was explained.
\bigskip

2. The Khintchine-Meinardus method used in this work covers a variety of models given by generating functions $\nicef(\delta)=f(e^{-\delta})$ exhibiting exponential asymptotics, as $\delta\to 0^+$ (see Lemma~\ref{estimates} below) which are essentially implied  by Meinardus Condition 1 of the main theorem. In this connection note that there exists  a rich literature (see e.g. \cite{flaj}) devoted to the case of moderately, i.e.~non exponentially, growing  $\nicef(\delta),\ \delta\to 0$. Such models are studied by a quite different singularity analysis.

\nin\section{Gentile statistics}
Gentile statistics is a model studied in physics
\cite{Gentile,SMB,TMB}, which
counts
partitions of an integer $n$ with no part occurring more
than $\eta-1$ times, where $\eta\geq 2$ is a parameter .
When $\eta=2$, Fermi-Dirac statistics are obtained and when $\eta=\infty$,
Bose-Einstein statistics, with uniform weights $b_k=1, \ k\ge 1$ result.
As far as we know,
no rigorous derivation of the asymptotics of Gentile statistics has
previously been given, although Theorem~\ref{gentile}
below was anticipated in approximation (23) of \cite{SMB}.
In this work we derive the aforementioned theorem as a special case of our Theorem~\ref{general}.

Gentile statistics are Taylor coefficients of the generating function
$$
f(z)=\prod_{k=1}^\infty\frac{1-z^{\eta k}}{1-z^k},\quad \vert z\vert<1,\quad \eta\ge 2 \ \,{\rm is \ an  \ integer}
.
$$
We remark that there is another natural interpretation of the Gentile
statistics, which is the number of integer partitions with no part
size divisible by $\eta$,
where part sizes can appear an unlimited number
of times.
Gentile statistics fit into the framework \refm[frame]
of Theorem~\ref{general} with
$$
S(z)=\frac{1-z^{\eta}}{1-z},\quad \vert z\vert<1,\quad \eta\ge 2
\ \,{\rm is \ an \ integer}
$$
and $a_k=b_k=1, \quad k\ge 1$.
\begin{theorem}\label{gentile}
Gentile statistics have asymptotics
$$
c_n\sim\sqrt{\frac{\kappa}{4\pi\eta}}n^{-3/4}e^{2\kappa\sqrt{n}},
$$
where
$$
\kappa=\sqrt{\zeta(2)(1-\eta^{-1})},\quad \eta\ge 2\ \,
{\rm is \ an \ integer}.
$$
\end{theorem}
\proof
We will show that  all the conditions  of Theorem~\ref{general} are satisfied
for Gentile statistics.
In order to show that
\refm[assum2]  holds for $\eta>1$, we calculate
$$
\frac{d^2}{d\delta^2}\log S(e^{-\delta})=
\frac{e^{\delta}}{(e^{\delta}-1)^2}
-\frac{\eta^2e^{\eta\delta}}{(e^{\eta\delta}-1)^2}.
$$
We have
\begin{eqnarray*}
\frac{d}{d\eta}\frac{\eta^2e^{\eta\delta}}{(e^{\eta\delta}-1)^2}
%&=&
%\frac{{(e^{\eta\delta}-1)^2}e^{\eta\delta}(2\eta+\delta\eta^2)
%-2\delta\eta^2e^{2\eta\delta}(e^{\eta\delta}-1)}{(e^{\eta\delta}-1)^4}\\
&=&
\frac{\eta e^{\eta\delta}\left[e^{\eta\delta}(2-\delta\eta)
-(2+\delta\eta)\right]}{(e^{\eta\delta}-1)^3}\\
&=&
\frac{\eta e^{\eta\delta}g(\eta\delta)}{(e^{\eta\delta}-1)^3},
\end{eqnarray*}
where $g(x)=e^x(2-x)-(2+x)$. Taking the derivative of $g$ produces
$g^\prime(x)=e^x(1-x)-1<0$ for $x>0$, which, together with $g(0)=0$,
implies that $g(x)<0$ for $x>0$. Combining this with the fact that $
\frac{d^2}{d\delta^2}\log S(e^{-\delta})=0,$ if $\eta=1,$ we conclude
that \refm[assum2] holds, for all $\eta>1.$

It remains to be shown that conditions $(I) - (III)$ are satisfied
for the model considered.
%We can make use of \refm[factor], because $a_k=1$.  We have $D_b(s)=\zeta(s)$, where $\zeta(s)$ is the Riemann zeta function.
We have
$$
\log f(z) = \sum_{k=1}^\infty\sum_{j=1}^\infty
\left(\frac{z^{jk}}{j}-\frac{z^{jk\eta}}{j}\right),\quad \vert z\vert<1,
$$
and so, by \refm[logflam], \refm[Dirdef] and \refm[dpres],
\begin{eqnarray*}
D(s)&=&\sum_{k=1}^\infty\sum_{j=1}^\infty
\left(\frac{(jk)^{-s}}{j}-\frac{(jk\eta)^{-s}}{j}\right)\\
&=&
\zeta(s)\zeta(s+1)(1-\eta^{-s}).
\end{eqnarray*}
Conditions $(I)$ and $(II)$ are satisfied
with any $0<C_0<1$
because of the analytic continuation of the Riemann
zeta function and the well known bound
\be\label{wellknown}
\zeta(x+iy)=O(|y|^C),\quad y\to \infty,
\en
for a constant $C>0,$ uniformly in $x$.
It is easy to check that $l_0=1$
and
$b_ka_k=1=k^{\rho_1-1}$, where $\rho_1=1$, and so \refm[assum1] is satisfied.
Hence condition $(III)$ is satisfied.
Moreover,
$$
r=1,\ \rho_0=0, \rho_1=1,\ A_0=\lim_{s\to 0} sD(s)= 0,\  A_1=\zeta(2)(1-\eta^{-1}),\\
$$
$$
\Theta=\lim_{s\to  0}D(s)= \zeta(0)\log\eta.
$$
By the argument preceding Proposition 1 in Appendix A, this says that the integrand $\delta_n^{-s}\Gamma(s)D(s)$ has a simple pole at $s=0$
 with residue $\Theta=\zeta(0)\log \eta$ and a simple pole at $s=1$ with residue $\zeta(2)(1-\eta^{-1})\delta_n^{-1}.$
As a result, in the case considered $\delta_n=\hat{h}_1^{1/2}n^{-1/2}-2^{-1}\hat{h}_0 n^{-1} + O(n^{-\frac{C_0}{2}-1})$,
where $\hat{h}_1$ and $\hat{h}_0$ are defined by 
(\ref{hhatldef}) and (\ref{hhat0def}),
and we arrive at the claimed asymptotic formula for $c_n$. \qed

\section{Asymptotic enumeration for distinct part sizes}

Weighted partitions fit our framework \refm[frame] with $S(z)=(1-z)^{-1}$,
$a_k=1,\ k\ge 1$ and weights $b_k$. When $b_k=1,\ k\ge 1$, Theorem~\ref{general} gives
the principal term in the asymptotical expansion of the number of partitions of $n,$ obtained by Hardy and Ramanujan in their famous paper \cite{HR}.
If
$S(z)=1+z$, $a_k=1,$ $b_k=k^{r-1},\ r>0,\ k\ge 1$,
then $c_n$ enumerates
weighted partitions having no repeated parts, called
expansive selections.
The asymptotics of expansive selections
were also studied in \cite{GS2}.

In this section, we find the asymptotics of $c_n$ induced by the generating function
\be\label{distinct}
f(z)=\sum_{n=0}^\infty c_n z^n= \prod_{k=1}^\infty(1+ k^{-q} z^k), \vert z\vert<1, \quad q>0.
\eu
The model fits the setting \refm[frame] with $S(z)=1+z,\ b_k=1,\  a_k=k^{-q},\ k\ge 1$ and it can be considered as a colored selection with parameter $k^{-q}$  proportional to the number $m_k$ of colors of a component of size $k,$ e.g. $m_k=  y^k k^{-q},$ for some $y>1.$
A particular case of the model, when $q=1$ was
studied in Section~ 4.1.6 of \cite{GK} where  it was proven, with the help of a Tauberian theorem, that in this case
\be \lim_{n\to\infty}c_n=e^{-\gamma} \la{q1}\eu
and it was established the rate of convergence of $c_n,\  n\to\infty.$
Also, in  \cite{GK} it was shown that $c_n$ is equal to the probability that a random polynomial  of order $n$ is a product of irreducible factors of different degrees.
In \cite{O}, Section 11, it was demonstrated  that $c_n$ can be treated as the probability that a random permutation on $n$ letters has distinct cycle   lengths, and another proof of \refm[q1] was suggested.
The generating function \refm[distinct]
is discussed in \cite{flaj} and the method therein
applied when $q>1$ to give the rate of decay $n^{-q}$ of $c_n$,
up to an unspecified constant.
By elementary techniques we show below
 that for $q>1,$ $c_n\sim W(q) n^{-q},$ where $W(q):=\sum_{n=0}^\infty c_n$.
In Appendix~\ref{Wexp}, an expression
for the constant $W(q)$ as an infinite  product is derived.

Finally, note that in \cite{O}, (11.35), it was shown that for $q=2$,
the generating function $f(z)$ can not be analytically continued beyond the unit circle.

\begin{theorem}\label{weighted}
Let
$$
\sum_{n=0}^\infty c_n z^n= \prod_{k=1}^\infty(1+ k^{-q} z^k),\ \vert z\vert<1.
$$
If $0<q<1$, then $c_n$ has asymptotics given by \refm[cnasymp] with
$r=\max\{j\ge 1: 1-qj>0\}$ and $\rho_l=1-ql,\ l=1,\ldots r.$
If $q>1$, then, for a constant $W(q)>0$ depending only on $q$,
\be\label{cnbounds}
c_n\sim W(q) n^{-q},
\eu
as $n\to\infty$.
\end{theorem}
\proof\\
\nin{\bf The case $0<q<1$}.\\
We will apply Theorem~\ref{general}.
Assumption \refm[assum2]
is easy to verify.
We have
$$
 \log f(z)=\sum_{k\ge 1} \log\left(1+\frac{z^k}{k^q}\right)=\sum_{k\ge 1} \sum_{j\ge 1}(-1)^{j-1} \frac{z^{kj}}{jk^{qj}},\quad \vert z\vert<1
$$
and so, by \refm[logflam] and \refm[dpres],
$$
D(s)=D(s;q)=\sum_{k\ge 1} \sum_{j\ge 1}(-1)^{j-1} \frac{{(kj)}^{-s}}{jk^{qj}}\\
=
\sum_{k\ge 1} \sum_{j\ge 1}
\frac{(-1)^{j-1}}{j^{s+1}k^{s+qj}}.
$$

We claim that the function $D(s;q)$ allows
analytic continuation to the set
$\BC$ except for the poles in $H_q:=\{s= 1- qj,\ j=1,2,\ldots,\ q<1 \}$.
Changing the order of summation, we write
\be
D(s;q)=\sum_{j\ge 1} \sum_{k\ge 1}\frac{(-1)^{j-1}}{j^{s+1}k^{s+qj}}=\sum_{j\ge 1} \frac{(-1)^{j-1}}{j^{s+1}}\zeta(s+qj),\quad \Re(s)> 0, s\notin H_q. \la{yurd}
\eu
Note that
$$\zeta(s+qj)=1+\sum_{n=2}^\infty \frac{1}{n^{s+qj}}:= 1 + \Phi(s;q),$$
where the function $\Phi(s;q)$ is analytic for $s\in \BC\setminus H_q ,$ and moreover $$\Phi(s;q)= O(2^{-qj}), \quad j\to \infty,\quad q>0,$$ uniformly in $s$ from any compact subset of $\BC\setminus H_q.$
This implies that the series
 $$\sum_{j\ge 1} \frac{(-1)^{j-1}}{j^{s+1}}\Phi(s;q)$$
converges absolutely and uniformly on any compact subset of $\BC\setminus H_q. $
By the Weierstrass convergence theorem, this implies
that the series above is analytic
in the above indicated domain. Since the function
$$\sum_{j\ge 1} \frac{(-1)^{j-1}}{j^{s+1}}= -(2^{-s}-1)\zeta(s+1)$$
is analytic in $\BC$, our claim is proven.
This allows us to conclude that condition $I$ of Theorem 1 holds with
$r=\max\{j\ge 1: 1-qj>0\}$ simple positive poles at
$\rho_l=1-q(r-l+1),\ l=1,\ldots r$ and with $0<C_0<1$ defined by
$$C_0=\left\{
        \begin{array}{ll}
          (r+1) q-1-\epsilon, \quad 0<\epsilon<(r+1)q-1, & \hbox{if } (r+1)q\le 2 \\
          {\rm any \ number \ in}\  (0,1), & \hbox{if } (r+1)q> 2.
        \end{array}
      \right.
$$

Condition $(II)$ follows from
\refm[wellknown] and \refm[yurd].
Finally,  $l_0=1$ in the case considered because $S(z)=1+z$ and
$b_ka_k=k^{-q}=k^{\rho_r-1}$ and so \refm[assum1] is satisfied.
Hence condition $(III)$ is satisfied, by Lemma 1 in \cite{GSE} which states that the bound \refm[assum1] on $b_k$ supplies
\refm[III].

The poles $\rho_l=1-ql$ are such that their differences are
multiples of $q$
and so $\Upsilon_r$ defined by \refm[sd2]
will contain the multiples of $q$ in $(0,\rho_r+1]$.
However, the set $\tilde\Upsilon_r$ defined by \refm[sd1]
may contain additional elements. For example, taking $\tilde d_0=2$
and all other $\tilde d_k=0$ in \refm[sd2] gives $2\rho_r=2(1-q)<2-q=\rho_r+1$
and so $2\rho_r\in\tilde\Upsilon_r$.
It will be the case that $3\rho_r\in\tilde\Upsilon_r$ if and only if $q>1/2$.
The set $\Upsilon_r$ defined by \refm[sd2]
will also contain $\rho_r$ (taking $k=0$).
Thus, the $\lambda_s$ will contain the multiples of $q$ and of $\rho_r=1-q$
in $(0,\rho_r+1]$. The coefficients $K_{s,l}$ in
\refm[cnasymp] are implicitly given by the recursion
$$
n\delta_n^{\rho_r+1}-\sum_{k=0}^r
\hat{h}_k\delta_n^{\rho_r-\rho_k}=o(n^{-1}),
$$
which comes from Proposition~1 of \cite{GS}. Note that 
 the definition of $\Upsilon_r$ is motivated by the solution of the above recursion.
To find the $K_{s,l}$ in examples such as this one, for which
$\delta_n$ has a series expansion,
one would probably want to
find the asymptotic solution of the recursion above by computer.

{\bf The case $q>1$.}\\
Theorem 1 is not applicable in this case, because  all  poles $1-qj,\ j\ge 1,\quad q>1$ of the function $D(s;q)$ in \refm[yurd], are negative.
From
\be f(z)=\prod_{k=1}^\infty (1+ z^kk^{-q})=\sum_{n=0}^\infty c_n z^n,\quad \vert z\vert\le 1,\ c_0=1 \quad q>1 \la{fprod}\eu
we have \be f(1)=\prod_{k=1}^\infty (1+ k^{-q}):=W(q)<\infty,\quad q>1,\la{jsa}\eu
since the convergence of the infinite product in \refm[jsa] is equivalent to the convergence of the series
$$\sum_{k=1}^\infty k^{-q}<\infty, \quad q>1.$$
By \refm[fprod] and \refm[jsa], \be W(q)=\sum_{n=0}^\infty c_n,\ q>1.\la{sumc}\eu

In proving that $\lim_{n\to \infty}n^qc_n=W(q)$, 
%we  use the common notation $$[z^k]f(z)=c_k$$ for the coefficient of $z^k$ of the function $f(z)$ in \refm[fprod] and denote by $$[z^k]f_l(z)=c_k^{(l)}$$ the same coefficient for the $l$-truncation $f_l$ of $f$ which is defined by
we define
$$
f_l(z):=\prod_{k=1}^l (1+ z^kk^{-q})
=\sum_{k=0}^{l(l+1)/2} c_k^{(l)}z^k.
$$
Two key facts  follow from the definitions of $f$ and $f_l$:
\be 0\le c_k^{(l)}\le c_k,\ 1\le l< k;\ \  c_k^{(l)}=c_k, \ l \ge k,\ k=0,1,\ldots l,
 \la{fact}\eu
Also note that $c_k^{(l)}$ counts the number of distinct partitions of $k$ with parts at most $l$ and that  \be  c_k^{(n-k-1)}=0,\ \hbox{if}\ \  k>\frac{(n-k-1)(n-k)}{2},\la{kcond}\eu because $f_l(z)$ is a polynomial in $z$ of degree
$\frac{l(l+1)}{2}.$
From \refm[fprod] we obtain the recurrence relation
\begin{align}
c_n= & n^{-q}c_0 + (n-1)^{-q}c_{1}+\cdots+ (n-n^*)^{-q}c_{n^*}\nonumber\\
& + (n-n^*-1)^{-q}c_{(n^*+1)}^{(n-n^*-2)}+\cdots+
  +2^{-q}c_{n-2}^{(1)}+ c^{(0)}_{n-1},\ \
   n\ge 2, \la{recur}
\end{align}
where, by virtue of the second fact in \refm[fact],
$$
n^*=
\left\{
  \begin{array}{l l}
    [\frac{n}{2}], & \hbox{if n is  odd;} \\

    [\frac{n}{2}] - 1, & \hbox{if n  is even.}
  \end{array}
\right.
$$
The condition \refm[kcond] is equivalent to $n-1\ge k>k_n^*:= n+\frac{1}{2} - \frac{1}{2}\sqrt{8n +1}$ and it expresses the fact that dictinct partitions of $k_n^*\sim n-\sqrt{2n}, \ n\to \infty$   with largest part at most $n-k_n^*\sim \sqrt{2n},\ n\to \infty$ do not exist.
Hence, \refm[recur] can be written as
\be n^{q}c_n= \sum_{k=0}^{n^{\epsilon}}\left(\frac{n-k}{n}\right)^{-q}c_k + \sum_{k=n^\epsilon+1}^{n^*}\left(\frac{n-k}{n}\right)^{-q}c_k+ \sum_{k=n^{*}+1}^{k_n^*}\left(\frac{n-k}{n}\right)^{-q}c_k^{(n-k-1)},\la{sums}\eu
for some $\epsilon>0$.
The expression \refm[sums] and the first fact of \refm[fact]
imply the bound
$$
n^qc_n\leq\left(\frac{n-k_n^*}{n}\right)^{-q}\sum_{k=0}^{k_n^*}c_k
\leq\left(\frac{n-k_n^*}{n}\right)^{-q}W(q)
=O(n^{q/2}), \quad n\to \infty,
$$
and so
$$
c_n=O(n^{-q/2}), \quad n\to \infty.
$$

%In the third sum of \refm[sums], $$\frac{1}{2}\sqrt{8n +1}@@-1< n-k< n-n^* ,$$
The first fact of \refm[fact] now produces $c_k^{(n-k-1)}\le c_k=O(k^{-q/2})$,
which implies that the third sum of \refm[sums] is of order
$$O(1)\sum_{k=n^{*}+1}^{k_n^*}\left(\frac{n-k}{n}\sqrt{k}\right)^{-q}c_k
=O(1)\Big(\frac{(n-k_n^*)\sqrt{k_n^*}}{n}\Big)^{-q} \sum_{k=n^{*}+1}^{k_n^*} c_k\to 0, $$
where the last step holds
because $\lim_{n\to\infty}\frac{(n-k_n^*)\sqrt{k_n^*}}{n}=\sqrt{2}$
and because the last  sum is  the tail of a convergent series. Similarly, in the second sum,
$\frac{n-k}{n}\geq\frac{n-n^*}{n}>1/2$, which implies that the
second sum also vanishes as $n\to \infty.$
As a result, the main contribution comes from the first sum, for which
$\frac{n-k}{n}\uparrow 1$ termwise as $n\to \infty$, so that
the limit of the first sum as $n\to\infty$ is
 $W(q)$ by dominated convergence.\qed

{\bf Remark:} Comparing the asymptotics of $c_n$ in the cases $0<q<1$,\ $q=1$ and $q>1$ it is clearly seen that $q=1$ is a point of phase transition.
%\nin {\bf Computation of $V(q)$}.

\appendixtitleon
\appendixtitletocon
\begin{appendices}
\section{The proof of Theorem~\ref{general}}\label{proofs}
The first step in the proof is to find  the asymptotics of
$\nicef(\delta):=f(e^{-\delta})$, as $\delta\to 0^+$, because that
will help us estimate $f_n(e^{-\delta})$ for an appropriate choice
of $\delta=\delta_n$.
This is done in Lemma~\ref{estimates} below by
using the Mellin transform method of Meinardus and Condition $(I)$
on the poles of $D(s)$.
In the probabilistic approach initiated by
Khintchine,
the free parameter $\delta=\delta_n$ is chosen to be the solution
of the equation
\be
\BE U_n=n,\ n\ge 1. \la{choice}
\eu

We introduce the notations
\be\label{hhatldef}
\hat{h}_l=\rho_lh_l, \quad l=1,\ldots r,
\eu
\be\label{hhat0def}
\hat{h}_0=-A_0,
\eu
where the $h_l$ are defined in \refm[hldef].
In Proposition~\ref{deltafacts} below
we use the fact that
the equation \refm[choice] for $\delta_n$ can be written as
\be\Big(-\log(f_n(e^{-\delta}))\Big)^\prime_{\delta=\delta_n}=n,\quad n\ge 1\la{der} \eu
to derive the facts  \refm[sd6], \refm[fnasymp], and \refm[bog].
In \cite{GS}, by a careful analysis of the equation \refm[der]  the following  expansion of  $\delta_n, \ n\to \infty$ was derived:
\be
\delta_n=\big(\hat{ h}_r\big)^{\frac{1}{\rho_r+1}}n^{-\frac{1}{\rho_r+1}}+ \sum_{s=1}^{\vert\Upsilon_r\vert}\hat{ K}_s n^{-\frac{1+\lambda_s}{\rho_r+1}}+o(n^{-1}), \la{sd6}
\en
where  $\hat{K}_s$ do not depend on $n$, and the powers $\lambda_s$ are as defined in \refm[sd2].
In \cite{GS} it was  also shown that
\be\label{fnasymp}
f_n(e^{-\delta_n})\sim \nicef(\delta_n),\quad n\to \infty,
\eu
and, moreover, that
\be
\left(\log f_n(e^{-\delta})\right)^{(k)}_{\delta=\delta_n}=
\left(\log \nicef(\delta)\right)^{(k)}_{\delta=\delta_n} +
\epsilon_k(n), \la{bog}\en
for $k=0,1,2,3$, where $\epsilon_k(n)\to 0,\quad n\to \infty,\quad k=0,1,2,3$.

We now analyze the three factors in the representation \refm[rep] when
$\delta=\delta_n$ is given by \refm[sd6].

$(i)$ It follows from \refm[sd6] that the first factor of \refm[rep]
equals
\be
e^{n\delta_n}=\exp\left\{\big(\hat{h}_r\big)^{\frac{1}{\rho_r+1}}n^{\frac{\rho_r}{\rho_r+1}}+ \sum_{s:\lambda_s\le \rho_r} \tilde{K}_s n^{\frac{\rho_r-\lambda_s}{\rho_r+1}}+ \epsilon_n\right\}, \la{ndel}
\en
where $\lambda_s\in \Upsilon_r$ and $ \epsilon_n\to 0,\ n\to\infty$.

$(ii)$
For $l=0,1,\ldots, r$
$$
\big(\delta_n\big)^{-\rho_l}= \big(\hat{h}_r\big)^{\frac{-\rho_l}{\rho_r+1}}n^{\frac{\rho_l}{\rho_r+1}} + \sum_{s:\lambda_s\le \rho_l} K_{s,l} n^{\frac{\rho_l-\lambda_s}{\rho_r+1}}+\epsilon_{n}(l),
$$
where
$\epsilon_{n}(l)\to 0,\quad n\to \infty,\ l=1,2,\ldots r$,
and where the coefficients $K_{s,l}$ are obtained from the binomial expansion for $\big(\delta_n\big)^{-\rho_l}$, based on \refm[sd6] and the definition \refm[sd2] of the set $\Upsilon_r$.
Consequently, substituting $\delta=\delta_n$ into the expression
\refm[prod1] for ${\cal F}(\delta)$ in Lemma~\ref{estimates} below gives
$$
\log f_n(e^{-\delta_n})= \sum_{l=0}^{r}\hat{h}_l\big(\hat{h}_r\big)^{\frac{-\rho_l}{\rho_r+1}}n^{\frac{\rho_l}{\rho_r+1}}+
\sum_{l=0}^r \hat{h}_l\sum_{s:\lambda_s\le \rho_l} K_{s,l} n^{\frac{\rho_l-\lambda_s}{\rho_r+1}}+$$
\be \Big(\frac{A_0}{\rho_r+1}\log n-  \frac{A_0}{\rho_r+1}\log\hat{h}_r\Big)+\epsilon_n,\ \epsilon_n\to 0,\ n\to \infty.
\la{logfn}
\en

$(iii)$
The asymptotics of the third factor of \refm[rep]
are given by Theorem~\ref{LLT},
which is a local limit theorem, and which is proved using
condition $(III)$. The proof uses ideas from \cite{frgr}.
As a result of Theorem~\ref{LLT}, we have
$$
\PP\left(U_n=n\right)
 \sim\frac{1}{\sqrt{2\pi
K_2}}\big(\hat{h}_r\big)^{\frac{2+\rho_r}{2(\rho_r+1)}}n^{-\frac{2+\rho_r}{2(\rho_r+1)}},\quad
 n\rightarrow\infty,\non
$$
for a constant $K_2=h_r\rho_r(\rho_r+1)$.

Finally, to completely account for the influence of all $r+1$  poles $\rho_0,\rho_{1},\ldots,\rho_r$, we present the sum of the expressions \refm[ndel],
\refm[logfn] obtained in (i),(ii) for the first two  factors in the representation \refm[rep]
in the following form:
\begin{eqnarray*}
n\delta_n+ \log f_n(e^{-\delta_n})&=&
\sum_{l=0}^{r}P_ln^{\frac{\rho_l}{\rho_r+1}}+
\sum_{l=0}^r h_l\sum_{s:\lambda_s\le \rho_l} K_{s,l} n^{\frac{\rho_l-\lambda_s}{\rho_r+1}}\\
&&+ \Big(\frac{A_0}{\rho_r+1}\log n -  \frac{A_0}{\rho_r+1}\log\hat{h}_r\Big )+\epsilon_n,
\end{eqnarray*}
where $P_l$ denotes the resulting coefficient of $n^{\frac{\rho_l}{\rho_r+1}}.$

If $r=1$, then \refm[der], \refm[exactly] produce
$$
n=\hat{h}_1\delta_n^{-\rho_r-1}+\hat{h}_0\delta_n^{-1}+O(\delta_n^{C_0-1})
+\varepsilon(n),
$$
with $\varepsilon(n)\to 0$, $n\to\infty$, which is analogous to
equation (54) of \cite{GSE}. The previous equation can be inverted
as in \cite{GSE}, giving
\be\label{lambdaasymp}
\delta_n=\hat{h}_1^{\frac{1}{\rho_1+1}}n^{-\frac{1}{\rho_1+1}}
+\frac{\hat{h}_0}{\rho_1+1}n^{-1}+O(n^{-1-\beta}),
\eu
where
$$
\beta=
\left\{
  \begin{array}{l l}
    \frac{C_0}{\rho_1+1} , & {\rm if}\ \ \rho_1\ge C_0; \\
    \frac{\rho_1}{\rho_1+1} , & {\rm otherwise}.
  \end{array}
\right.
$$
Substituting \refm[lambdaasymp] into the previous asymptotic estimates of
the three factors in \refm[rep], obtained in $(i)-(iii)$,
proves the asymptotic formula \refm[cnasymp] for $c_n$.
\bigskip

We now present the results to which we refered in the summary above.

\begin{lemma} \label{estimates}
\noindent  (i) As $\delta\to 0^+$,
\be \label{prod1}
\nicef(\delta)=\exp\left(\sum_{l=0}^r h_l\delta^{-\rho_l}
-A_0\log\delta+ M(\delta,C_0)\right),\en
%\end{document}
where $0<C_0<1,$
$\rho_0=0$,
\begin{eqnarray}
h_l&=&A_l\Gamma(\rho_l), \quad l=1,\ldots,r,\label{hldef}\\
h_0&=&\Theta-\gamma A_0,\nonumber\\
M(\delta;C_0)&=&\frac{1}{2\pi
i}\int_{-C_0-i\infty}^{-C_0+i\infty} \delta^{-s}\Gamma(s)D(s)ds=O(\delta^{C_0}),\ \delta\to 0,\nonumber
\end{eqnarray}
and where $\Theta$ is as in \refm[THETA].

(ii) The  asymptotic expressions for the derivatives
$$\Big(\log\nicef(\delta)\Big)^{(k)}$$
are given by the formal differentiation of the logarithm of
\refm[prod1], with \\  $(M(\delta;C_0))^{(k)}_\delta=O(\delta^{C_0-k}),\ k=1,2,3,\ \delta\to 0.$
\end{lemma}
\proof
We use the fact that
$e^{-u}$, $u>0$, is the Mellin transform of the Gamma function:
\be e^{-u}=\frac{1}{2\pi i}\int_{v-i\infty}^{v+i\infty}
u^{-s}\Gamma(s)\,ds,\quad u>0,\ \Re(s)=v>0. \la{Mellin} \en
Applying \refm[Mellin] with $v=\rho_r+\epsilon,\ \epsilon>0$ we have\begin{eqnarray}
\log
~\nicef(\delta)
&=&
\sum_{k=1}^\infty b_k\log S\left(a_k e^{-\delta k}\right)\nonumber\\
&=&
\sum_{k=1}^\infty b_k\sum_{j=1}^\infty \xi_j a_k^j e^{-\delta jk}\nonumber\\
&=&
\sum_{k=1}^\infty\sum_{j=1}^\infty b_k \xi_j a_k^j
\frac{1}{2\pi i}\int_{\epsilon+\rho_r-i\infty}^{\epsilon+\rho_r+i\infty}
(\delta jk)^{-s}\Gamma(s)\,ds\nonumber\\
&=&
\frac{1}{2\pi i}\int_{\epsilon+\rho_r-i\infty}^{\epsilon+\rho_r+\infty}
\delta^{-s}
\Gamma(s)
\sum_{k=1}^\infty\sum_{j=1}^\infty b_k \xi_j a_k^j
(jk)^{-s}\,ds\nonumber\\
&=&
\frac{1}{2\pi i}\int_{\epsilon+\rho_r-i\infty}^{\epsilon+\rho_r+i\infty} \delta^{-s}\Gamma(s)D(s)ds, \la{intrep3}
\end{eqnarray}
where we have used \refm[logflam],\refm[Dirdef] and \refm[logSexp]-\refm[dpres] at \refm[intrep3].
Next, we apply the residue theorem for the integral \refm[intrep3], in the complex domain $ -C_0\le \Re( s)\le \rho_r+\epsilon, $ with 
$0< C_0<1, \epsilon>0.$
By virtue of  condition $(I)$,
 the integrand in \refm[intrep3] has in the above domain $r$
simple poles at $\rho_l>0, \ l=1,\ldots,r$.
The corresponding residues at $s=\rho_l$ are equal to:
$A_l\delta^{-\rho_l}\Gamma(\rho_l),$ $l=1,\ldots,r$.

By the Laurent expansions at $s=0$ of the
Gamma function $\Gamma(s)=\frac{1}{s}-\gamma+\ldots,$
and the function $D(s)=\frac{A_0}{s}+\Theta+\cdots$,
 the integrand $ \delta^{-s}D(s)\Gamma(s)$ may also have
a pole at $s=0$, which is a simple one with residue $\Theta,$ if $A_0=0,\Theta\neq 0,$
and is of  a second order with residue $\Theta-\gamma A_0 - A_0\log\delta$,
if $A_0\neq 0$.
In the case $A_0=\Theta=0$, the integrand $\delta^{-s}D(s)\Gamma(s)$ is analytic at $s=0$.

To apply the residue theorem, we bound the aforementioned domain by two horizontal contours $\vert \Im(s)\vert=t>0.$
 By condition $(II),$ the integral of the integrand $\delta^{-s}\Gamma(s)D(s)$, over the horizontal contours $-C_0\le\Re(s)\le \epsilon+\rho_r$, $\vert\Im(s)\vert=t>0$, tends to zero, as $t\to \infty,$ for
any fixed $\delta>0.$   This gives the claimed formulae \refm[prod1], where the remainder term $M(\delta;C_0)$ is the
integral taken over the vertical contour $-C_0+it,\ -\infty<t<\infty$.
This proves $(i)$.

In order to prove $(ii)$,
%$\left(\log\nicef(\delta)\right)^{(1)}$,
one differentiates
the logarithm of \refm[prod1] with respect to $\delta$ and estimates
the remaining integral in the same way as above.
\hfill\qed

We will need the following bound on $b_k$.
\begin{proposition}
Let the double series
%\be D(s)=\sum_{k,j=1}^\infty b_k \xi_j a_k^j (jk)^{-s} \la{double}\eu
$D(s)$ defined by \refm[dpres]
converge absolutely in the half-plane ${\cal R}(s)>\rho,$ for some $\rho>\rho_r.$
Then
the following bound holds \be b_ka_k^{l_0}=o(k^{\rho}),\quad k\to \infty,\la{bound2}\eu
where $l_0$ is defined by \refm[l0def].
\end{proposition}
\proof
We observe that $l_0$ defined by \refm[l0def] satisfies
\be
l_0=\min\{j\ge 1:\xi_j\neq 0\}.\la{l0}\eu

The assumed absolute convergence of the double series in \refm[dpres]
implies the  absolute convergence for $\rho>\rho_r$
of the  iterated series
$$\sum_{j=1}^\infty \frac{\xi_j }{j^\rho}\sum_{k=1}^\infty\frac{b_ka_k^j}{k^\rho}.$$
Consequently,$$\sum_{k=1}^\infty\frac{b_ka_k^j}{k^\rho}<\infty,\quad
{\rm for \ all}\quad  j\ge 1: \xi_j\neq 0.$$
Hence,$$\frac{b_ka_k^j}{k^\rho}\to 0,\quad k\to \infty,\quad {\rm for \ all}\quad  j\ge 1: \xi_j\neq 0.$$
The latter implies \refm[bound2]. \qed

\begin{proposition}\label{deltafacts}
For sufficiently large $n$ there is a unique solution of \refm[choice]
and that solution satisfies \refm[sd6], \refm[fnasymp], and \refm[bog].
\end{proposition}
\proof
For each $n\ge 1$, the function
$\Big(-\log(f_n(e^{-\delta}))\Big)^\prime_\delta$
 is decreasing for all $\delta>0$ because
of \refm[assum2].
Moreover, setting
\be\label{deltadef}
\delta=\delta(n):=Cn^{-\frac{1}{\rho_r+1}}, C>0,
\eu
 we
have
\begin{eqnarray}
\Big(-\log(f_n(e^{-\delta}))\Big)^\prime_{\delta(n)}
&=&(-\log\nicef(\delta))^\prime_{\delta(n)}- \sum_{k=n+1}^\infty
\Big(- b_k\log S(a_ke^{-k\delta})\Big)_{\delta(n)}^\prime\nonumber\\
&=&(-\log\nicef(\delta))^\prime_{\delta(n)}- O\Big(
\sum_{k=n+1}^\infty b_k \xi_{l_0}a_k^{l_0}e^{-\delta kl_0}kl_0\Big)\nonumber\\
&\sim&C^{\frac{1}{\rho_r+1}}\rho_rh_r n,\quad C>0.
\la{deltas}
\end{eqnarray}
Here the step before the last follows because for $\delta$ defined by
\refm[deltadef] we have $n\delta=Cn^{\frac{\rho_r}{\rho_r+1}}\to \infty,\ n\to \infty,$  because of Lemma~\ref{estimates} $(ii)$   and due to the fact that
$$-\Big(\log S(a_ke^{-k\delta})\Big)_{\delta(n)}^\prime\sim \xi_{l_0}a_k^{l_0}e^{-\delta(n) kl_0}kl_0,\quad  n\to \infty,$$ uniformly for $k\geq n+1$,
where $l_0$ satisfies \refm[l0]. The last step in \refm[deltas] results from  Lemma~\ref{estimates} $(ii)$ and  \refm[bound2].
For $n$ sufficiently large, the right hand side of \refm[deltas] is $>n,$ if $C>(\rho_rh_r)^{-(\rho_r+1)}$ and $\le n$ otherwise.
This and \refm[assum2] say that for a sufficiently large $n,$ \refm[der] has a unique solution $\delta_n,$ which satisfies \be \delta_n\sim (\rho_rh_r)^{-(\rho_r+1)}n^{-\frac{1}{\rho_r+1}},\quad n\to\infty, \la{delasy} \eu
where $h_r$ is defined as in \refm[prod1].
We proceed to find an asymptotic expansion for $\delta_n$ by using a refinement of  the scheme of
Proposition~1 of \cite{GS}.
We call any $\tilde{\delta}_n,$ such that \be \big(-\log f_n(e^{-\delta})\big)_{\delta=\tilde{\delta}_n}^\prime -n\to 0,\quad n\to \infty \la{asympt}\en
an asymptotic solution of \refm[der].
We will show that it is sufficient for
\refm[asympt] that $\tilde{\delta}_n$ obeys the condition
\be (- \log \nicef(\delta))_{\delta=\tilde{\delta}_n}^\prime -n\to 0, \ n\to \infty. \la{aut} \en
By Lemma~\ref{estimates}, we have
$$ \big(- \log \nicef(\delta)\big)^\prime_\delta\sim h_r\rho_r \delta^{-\rho_r-1},\ \delta\to 0,$$
so that  \refm[aut] implies
\be\tilde{\delta}_n\sim (h_r\rho_r)^{\frac{1}{\rho_r+1}}n^{-\frac{1}{\rho_r+1}}, \ n\to \infty.\la{caf}\en
Next we have for all $n\ge 1$
\be
\log f_n(e^{-\tilde{\delta}_n})=\log \nicef(\tilde{\delta}_n)-\sum_{k=n+1}^\infty
b_k\log S(a_k e^{-k\tilde{\delta}_n}).\la{ser}\en
Applying the same argument as in \refm[deltas],
we derive the bound
\be
\sum_{k=n+1}^\infty
\Big(- b_k\log S(a_ke^{-k\delta})\Big)_{\delta=\tilde{\delta}_n}^\prime=
  o(1),\quad n\to \infty.
\la{arg}
\eu
Now, \refm[ser] and \refm[arg] show that \refm[aut] implies
\refm[asympt].
We will now demonstrate
that the error of   approximating the exact solution $\delta_n$ by the asymptotic solution $\tilde{\delta}_n$ is of order $o(n^{-1}).$
By the definitions of $\delta_n,$ $\tilde{\delta}_n$
 we have
\be \Big(-\log f_n(e^{-\delta})\Big)^\prime_{\delta=\delta_n}- \Big(-\log f_n(e^{-\delta})\Big)^\prime_{\delta=\tilde{\delta}_n}=\epsilon_n,\quad \epsilon_n\to 0,\quad n\to \infty.\la{epsli}\eu
Next, applying the Mean Value Theorem, we obtain
\be \left| \Big(-\log f_n\big(e^{-\delta_n}\big)\Big)^\prime-\Big(-\log f_n\big(e^{-\tilde{\delta}_n}\big)\Big)^\prime\right|=
\left|(\delta_n-\tilde{\delta}_n) \Big(\log f_n(e^{-u_n})\Big)^{\prime\prime}\right|,\la{ghj}\en
where
$$u_n\in[\min(\delta_n,\tilde{\delta}_n),\max(\delta_n,\tilde{\delta}_n)].$$ By  \refm[epsli],  the left hand side of \refm[ghj] tends to $0$, as $n\to \infty,$ while, by virtue of \refm[delasy],\refm[caf],
\be\Big(\log f_n(e^{-\delta})\Big)^{\prime\prime}_{u_n}\sim \rho_r(\rho_r+1)h_r(\delta_n)^{-\rho_r-2}= O(n^{\frac{\rho_r+2}{\rho_r+1}}),\quad n\to \infty.\la{second}\eu
Combining \refm[ghj] with \refm[second], gives the desired estimate
\be \left|\delta_n-\tilde{\delta}_n\right|=o(n^{-1}),\quad n\to \infty. \la{nmi}\en
An obvious modification of the argument in \refm[deltas] allows also to conclude that  $$\sum_{k=n+1}^\infty
b_k\log S(a_k e^{-k\delta_n})\to 0, \quad n\to \infty.$$
As a result, \refm[fnasymp] is valid.

By $(ii)$ of Lemma~\ref{estimates} we have
\be\label{exactly}
\Big(-\log\nicef(\delta)\Big)^{\prime}_{\delta}= \sum_{l=0}^r \hat{h}_l\delta^{-\rho_l-1}+ \Big(M(\delta;C_0)\Big)^\prime_{\delta}.
\en
This is exactly the starting point of the analysis of $\tilde{\delta}_n$
in Proposition~1 of \cite{GS}.
We may therefore apply Proposition~1 of \cite{GS}, which provides
detailed asymptotics for $\delta_n$, and \refm[nmi], to  conclude
that \refm[sd6] holds.

Finally, by an argument similar to the one for the proof of
\refm[arg] we see that \refm[bog] holds, as well.
\hfill\qed

The following estimate is central to our proof of Theorem~\ref{LLT}.
\begin{proposition}\label{useful}
Recall that $\phi_n(\alpha)$ is defined by \refm[phidef] and that $l_0$ is defined by
\refm[l0def]. Then
we have for all $\alpha\in {\cal R}$,
\begin{eqnarray}
\log|\phi_n(\alpha)|&=&\log|\phi_n(\alpha;\delta_n)|=
-2\sum_{k=1}^{\infty}\Lambda_k e^{-k\delta_n}\sin^2(\pi k\alpha)+
\epsilon_n\la{phieq}\\
&\leq&
-\frac{2d_{l_0}}{S^2(e^{-1/8l_0})}\sum_{k=(8l_0\delta_n)^{-1}}^n b_ka_k^{l_0}
e^{-\delta_n l_0 k}\sin^2(\pi \alpha l_0 k),
\la{phiupper}
\end{eqnarray}
where $\epsilon_n\to 0$ as $n\to\infty$.
\end{proposition}
\proof
We write $\log|\phi_n(\alpha)|,\ \alpha\in {\cal R}$, as
\begin{eqnarray}
\log|\phi_n(\alpha)|&=&
\frac{1}{2}
\sum_{k=1}^n b_k
\left\{
\log S(a_k e^{2\pi ik\alpha-k\delta_n})+\log S(a_k e^{-2\pi ik\alpha-k\delta_n})
-2\log S(a_k e^{-k\delta_n})\right\}\nonumber\\
&=&
\frac{1}{2}
\sum_{k=1}^\infty b_k
\left\{
\log S(a_k e^{2\pi ik\alpha-k\delta_n})+\log S(a_k e^{-2\pi ik\alpha-k\delta_n})
-2\log S(a_k e^{-k\delta_n})\right\}\nonumber\\
&&+O\left(\sum_{k=n+1}^\infty b_k\xi_{l_0}a_k^{l_0}e^{-k\delta_nl_0}\right)
\la{first}\\
&=&
\frac{1}{2}
\sum_{k=1}^\infty b_k
\left\{
\log S(a_k e^{2\pi ik\alpha-k\delta_n})+\log S(a_k e^{-2\pi ik\alpha-k\delta_n})
-2\log S(a_k e^{-k\delta_n})\right\}\nonumber\\
&&+o\left(\sum_{k=n+1}^\infty k^\rho e^{-k\delta_nl_0}\right)\label{reas1}\\
&=&
\frac{1}{2}
\sum_{k=1}^\infty b_k
\sum_{j=1}^\infty\xi_j a_k^j e^{-jk\delta_n}
\left(e^{2\pi ijk\alpha}+e^{-2\pi ijk\alpha}-2\right)
+ \epsilon_n
\label{reas2}\\
&=&
-2\sum_{k=1}^\infty
\sum_{j=1}^\infty b_k\xi_j a_k^j e^{-jk\delta_n}
\sin^2(\pi jk\alpha)+ \epsilon_n \nonumber\\
&=&
-2\sum_{\ell=1}^\infty\Lambda_\ell e^{-\ell\delta_n}\sin^2(\pi \ell\alpha)+ \epsilon_n\label{reas},
\end{eqnarray}
where $\epsilon_n\to 0,\ n\to \infty,$ \refm[first] and \refm[reas1] use
\refm[bound2],
 \refm[reas2] uses \refm[delasy] which implies
$k^\rho e^{-k\delta_nl_0}\to 0,\ n\to \infty,\ k\ge n+1 \hbox{\ and any}\ \rho>0.$ Finally, \refm[reas] follows from \refm[logflam].

As for the  inequality \refm[phiupper],
defining $\tau$ to be $\tau=\delta_n-2\pi i \alpha,\ \alpha\in\BR$, we have
\begin{eqnarray}
\log|\phi_n(\alpha)|&=&
\Re\left(\log f_n(e^{-\tau})-
\log
f_n(e^{-\delta_n})\right)\nonumber\\&=&
\frac{1}{2}\sum_{k=1}^n b_k \log\frac{\vert S(a_ke^{-k\tau})\vert^2}{S^2(a_ke^{-k\delta_n})}\nonumber\\ &=&
\frac{1}{2}\sum_{k=1}^n b_k \log\bigg(1-\frac{S^2(a_ke^{-k\delta_n})-\vert S(a_ke^{-k\tau})\vert^2}{S^2(a_ke^{-k\delta_n})}\bigg)\nonumber\\ &\le&
-\frac{1}{2}\sum_{k=1}^n b_k \frac{S^2(a_ke^{-k\delta_n})-\vert S(a_ke^{-k\tau})\vert^2}{S^2(a_ke^{-k\delta_n})},\la{vnbeg}
\end{eqnarray}
where the last inequality is because $S^2(a_ke^{-k\delta_n})-\vert S(a_ke^{-k\tau})\vert^2\ge 0,$ for all $\alpha\in {\cal R}$
and because $\log(1-x)\le -x,\ 0<x<1.$
Recalling \refm[sz] and \refm[l0def], we obtain for all $\alpha\in{\cal R},$\begin{eqnarray*}
S^2(a_ke^{-k\delta_n})-\vert S(a_ke^{-k\tau})\vert^2&=&
4\sum_{0\leq l,m<\infty}
d_ld_ma_k^{l+m} e^{-(l+m)k\delta_n} \sin^2\left((l-m)\pi\alpha k\right)\\
&\geq&
4d_{l_0}a_k^{l_0}e^{-\delta_n l_0k}\sin^2(\pi \alpha l_0k),\ k=1,2,\ldots,
\end{eqnarray*}
which allows us to continue \refm[vnbeg], arriving at the desired bound:
\begin{eqnarray*}
 \log|\phi_n(\alpha)|&\le&-2d_{l_0}\sum_{k=1}^n b_k a_k^{l_0}\frac{e^{-\delta_nl_0 k}\sin^2(\pi\alpha l_0k)}{S^2(a_ke^{-k\delta_n})}\\
&\le&
-2d_{l_0}\sum_{k=(8l_0\delta_n)^{-1}}^n b_ka_k^{l_0} \frac{e^{-\delta_n l_0k}\sin^2(\pi \alpha l_0k)}{S^2(a_ke^{-k\delta_n})}\\
&\le&
-\frac{2d_{l_0}}{S^2(e^{-1/8l_0})}\sum_{k=(8l_0\delta_n)^{-1}}^n b_ka_k^{l_0}
e^{-\delta_n l_0 k}\sin^2(\pi \alpha l_0 k),
\end{eqnarray*}
where the last step is because  $d_l\ge 0, \ l=1,2,\ldots,$  because $0<a_k\le 1$ and because $1\le S(z)<\infty$ is monotonically increasing  in $ \ 0\le z<1.$
\hfill\qed
\bigskip

\begin{theorem}[Local Limit Theorem]\label{LLT}

Let the random variable $U_n$ be defined as in \refm[fndef], \refm[phidef]. Then

\begin{eqnarray}
\PP\left(U_n=n\right) &
 \sim&\frac{1}{\sqrt{2\pi {\rm Var(U_n)
 }}}\non\\
&\sim&\frac{1}{\sqrt{2\pi
K_2}}\left(\delta_n\right)^{1+\rho_r/2}\non\\
& \sim&\frac{1}{\sqrt{2\pi
K_2}}\big(\hat{h}_r\big)^{\frac{2+\rho_r}{2(\rho_r+1)}}n^{-\frac{2+\rho_r}{2(\rho_r+1)}},\quad
 n\rightarrow\infty,\non
\end{eqnarray}
for a constant $K_2=h_r\rho_r(\rho_r+1)$.
\end{theorem}

{\bf Proof}\
\ We take $\delta=\delta_n$
in \refm[phidef] and
define
\be \alpha_0=\alpha_0(n)=(\delta_n)^{\frac{\rho_r+2}{2}}\log n.\la{alpha0}\eu
We write
$$
 \PP(U_n=n)=\int_{-1/2}^{1/2}\phi_n(\alpha)
e^{-2\pi in\alpha}d\alpha=I_1+I_2,$$
 where
$$
I_1=\int_{-\alpha_0}^{\alpha_0}\phi_n(\alpha)e^{-2\pi
in\alpha}d\alpha
$$
and
$$
I_2=\int_{-1/2}^{-\alpha_0}\phi_n(\alpha)e^{-2\pi
in\alpha}d\alpha +\int_{\alpha_0}^{1/2}\phi_n(\alpha)e^{-2\pi
in\alpha}d\alpha.
$$
The proof has two parts corresponding to  evaluation of the integrals $I_1$ and $I_2$, as $n\to \infty.$

{\bf Part 1:} Integral $I_1.$\ \
Defining $B_n$ and $T_n$ by \be
  B_n^2=
\Big(\log f_n\big(e^{-\delta}\big)\Big)^{\prime\prime}_{\delta=\delta_n} \ {\rm and}
\quad T_n=-\Big(\log f_n\big(e^{-\delta}\big)\Big)^{\prime\prime\prime}_{\delta=\delta_n}\label{BT2}
\end{equation}
for $n$ fixed we have the expansion in $\alpha$
%$$
%\varphi_j(t)=1+\BE{Y_j}t-
%\BE{Y_j^2}\frac{t^2}{2}+O\left(\BE{Y_j^3}t^3\right),\quad t\rightarrow0\\
%$$
 \begin{eqnarray}
 \non\phi_n(\alpha)e^{-2\pi in\alpha}
%&=&
%\exp{\left(\sum_{j=1}^{n}
%{\log{\left(1+\E{Y_j}t-\frac{t^2}{2}\BE{Y_j^2}+
%\BE{Y_j^3}O(t^3)\right)} -e^{-int}}\right)}\nonumber\\
&=&
\exp{\left(2\pi i\alpha(\BE U_n-n)-2\pi^2\alpha^2B_n^2+O(\alpha^3)T_n\right)}\nonumber\\
&=&\exp{\left(-2\pi^2\alpha^2B_n^2+O(\alpha^3) T_n\right)},\quad
\alpha\rightarrow0, \non
 \end{eqnarray}
where the second equation is due to \refm[choice].   By virtue of \refm[prod1] and \refm[bog] we derive from \refm[BT2] that  the main terms in the asymptotics for $B_n^2$ and   $T_n$ depend on the rightmost pole $\rho_r$ only:
\be B_n^2\sim K_2(\delta_n)^{-\rho_r-2}, \la{basym}\eu
$$T_n\sim K_3(\delta_n)^{-\rho_r-3},\quad n\to \infty$$
where $K_2=h_r\rho_r(\rho_r+1)$  and
$K_3=h_r\rho_r(\rho_r+1)(\rho_r+2)$ are  obtained from \refm[BT2] and Lemma 1.
Therefore, by virtue of \refm[alpha0],$$B^2_n\alpha_0^2\to \infty, \ T_n\alpha_0^3\to 0, \ n\to \infty.$$
Consequently, in the same way as in the proof of the local limit theorem in  \cite{GSE},
\be\label{I1sim}
I_1\sim \frac{1}{\sqrt{2\pi B_n^2}}, \ n\to \infty,
\en
and it is left to show that \be I_2=o(I_1), \ n\to \infty.\la{relat}\eu

{\bf Part 2:} Integral $I_2.$\ \
We rewrite the upper bound in \refm[phiupper] in Proposition~\ref{useful} as
\begin{eqnarray*}
\log|\phi_n(\alpha)|\le
-CV_n(\alpha),\quad \alpha\in \BR,
\end{eqnarray*}
where $C>0$ does not depend on $n$ and
$$
V_n(\alpha):= \sum_{k=(8l_0\delta_n)^{-1}}^n b_ka_k^{l_0}
e^{-\delta_n l_0 k}\sin^2(\pi \alpha l_0 k).
$$
It is enough to consider $\alpha\ge 0$, as the case $\alpha<0$ is exactly the same.

We split the interval of integration $[\alpha_0,1/2]$ into subintervals:
$$
[\alpha_0,(2\pi)^{-1}\delta_n]\cup
[(2\pi) ^{-1}\delta_n,1/2]{\rm \ \ if \ }l_0=1
$$ and
$$
[\alpha_0,(2\pi l_0)^{-1}\delta_n]\cup
[(2\pi l_0)^{-1}\delta_n, (2l_0)^{-1}]\cup
[(2l_0)^{-1},1/2] {\rm \ \ if \ }l_0>1.
$$
Our goal is to bound, as $n\to \infty,$  the  function $V_n(\alpha)$ from below in each of the subintervals. Firstly, we show that on the first two subintervals for $l_0\ge 1,$ the desired bound is implied by the
assumption \refm[assum1] in condition $(III)$.

In the first subinterval
$[\alpha_0,(2l_0)^{-1}\delta_n],\ l_0\ge 1$
we will  use the inequality \be
\sin^2(\pi x) \geq 4\parallel x\parallel^2,\ x\in\BR,  \la{th}
\en
where $\parallel x\parallel$
denotes the distance from $x$ to the nearest integer, i.e.
$$
\parallel x\parallel=
\left\{
\begin{array}{l l}
\{x\}&{\rm if \ }\{x\}\leq 1/2;\\
1-\{x\}&{\rm if \ }\{x\}>1/2.
\end{array}
\right.
$$
(see \cite{frgr} for the proof of \refm[th]).
By \refm[assum1] and \refm[th], we then have
\be V_n(\alpha) \ge 4 e^{-1/2}\sum_{k=(4l_0\delta_n)^{-1}}^{(2l_0\delta_n)^{-1}}C_2 k^{\rho_r-1}\parallel \alpha l_0 k\parallel^2,\ \alpha\in\BR,\ l_0\ge 1.
\la{ael}
\en
In the first subinterval,
$$
\parallel \alpha l_0k\parallel=\alpha l_0k, \ 1\le k\le (2l_0\delta_n)^{-1},\ l_0\ge 1,
$$
so that
\refm[ael] and \refm[alpha0] produce
\begin{eqnarray*}
V_n(\alpha)&\ge&
 4 C_2 e^{-1/2}l_0^2\alpha_0^2\sum_{k=(4l_0\delta_n)^{-1}}^{(2l_0\delta_n)^{-1}}
k^{\rho_r+1}\\
&\sim&
4C_2e^{-1/2}(\rho_r+2)^{-1}l_0^2 \alpha_0^2((2l_0\delta_n)^{-\rho_r-2}-(4l_0\delta_n)^{-\rho_r-2}) ,\quad  n\to \infty.
\end{eqnarray*}
By \refm[alpha0] and \refm[phiupper] this gives the desired bound in the first subinterval:
\be \log|\phi_n(\alpha)|\le - C\log^2 n,\quad C>0,\quad n\to \infty.\la{bou1}\en

For the second subinterval we will apply the argument in the proof of Lemma~1 in \cite{GSE}.
Given $\alpha\in\BR \setminus \{0\}$,
define the function $P(\alpha,\delta_n)$ by
\be P(\alpha,\delta_n)=\left[\frac{1+\vert \alpha
\vert \delta_n^{-1}}{2\vert \alpha\vert}\right]\ge 1,\la{P}\eu where
$[x]$ denotes the integer part of  $x$ and the inequality holds for $n$ large
enough.  The lower bound  for the zeta sum (see (17) in \cite{GSE}), supplies the bound
%\be\label{bnd0} \sum_{k=1}^P\sin^2(\pi k\alpha)\ge
%\frac{\delta_n^{-1}}{4}, \eu provided
\be
\sum_{k=(8\delta_n)^{-1}}^{P(\alpha,\delta_n)}\sin^2(\pi k\alpha)\ge
\frac{\delta_n^{-1}}{8}, \la{sinbo}
\eu
provided
\be
\frac{\delta_n}{2\pi}\leq\alpha\leq 1/2.
\la{alboun}\eu

Observing, that by definition \refm[P], $n>P(l_0\alpha,\delta_n)> (2\delta_n)^{-1}>(8\delta_n)^{-1}$
we rewrite \refm[sinbo] as

\be
\sum_{k=(8l_0\delta_n)^{-1}}^{P(l_0\alpha,\delta_n)}\sin^2(\pi kl_0\alpha)\ge
\frac{\delta_n^{-1}}{8}, \ l_0\ge 1,\la{sinbo1}
\eu
for
\be
(2\pi l_0)^{-1}\delta_n\leq\alpha\leq (2l_0)^{-1},\quad l_0\ge 1.
\la{alboun1}\eu
We now write $P$ for $P(l_0\alpha,\delta_n)$.
For $\alpha\in\BR$ and $l_0\ge 1$, we have
\begin{eqnarray*}
\sum_{k=(8l_0\delta_n)^{-1}}^n b_ka_k^{l_0}
e^{-\delta_nl_0 k}\sin^2(\pi \alpha l_0 k)
&\geq&
\sum_{k=(8l_0\delta_n)^{-1}}^P C_2 k^{\rho_r-1} e^{-\delta_nl_0 k}\sin^2(\pi \alpha l_0 k)\\&\geq&
C_2 e^{-Pl_0\delta_n}\sum_{k=(8l_0\delta_n)^{-1}}^P k^{\rho_r-1}\sin^2(\pi \alpha l_0k)\\
&:=&Q(\alpha).
\end{eqnarray*}
In order to get the needed lower bound on $Q(\alpha)$, we take into account that for all $\alpha$ obeying \refm[alboun1], $\frac{1}{2}<P\delta_n<\frac{1}{2}(1+2\pi ):=d.$
Applying  \refm[sinbo1], we  distinguish between the following two
cases: $(i)\ 0<\rho_r<1$ and $(ii)\ \rho_r\ge 1$.

For $\alpha$ in \refm[alboun1], we have  in  case $(i),$ $$Q(\alpha)\ge C_2
e^{-Pl_0\delta_n}P^{\rho_r-1}(8\delta_n)^{-1}\ge \frac{C_2}{8}
e^{-dl_0}(P\delta_n)^{\rho_r-1}\delta_n^{-\rho_r}\ge \frac{C_2}{8} d^{\rho_r-1}
e^{-dl_0}\delta_n^{-\rho_r}:= C_3\delta_n^{-\rho_r},$$ and in  case $(ii),$
$$Q(\alpha)\ge
C_2
e^{-Pl_0\delta_n}(8\delta_n)^{-\rho_r+1}(8l_0\delta_n)^{-1}\geq
C_2e^{-dl_0}8^{-\rho_r}\delta_n^{-\rho_r}l_0^{-1}:=C_4\delta_n^{-\rho_r}.$$
Finally, combining this with \refm[phiupper] gives the desired upper bound on $\log|\phi_n(\alpha)|$ for all $\alpha$ in \refm[alboun1] and $n$ sufficiently large:
\be \log|\phi_n(\alpha)|\le - C\delta_n^{-\rho_r},\ C>0.\la{logbou}\eu

{\bf Remark:}
If $l_0>1$, then $\sin^2(\pi \alpha l_0 k)=0, \ k\ge 1$ when $\alpha=l_0^{-1}\le 1/2,$
so that in  the third subinterval $[(2l_0)^{-1},1/2]$ the above bounds are not applicable.

In the third  subinterval $[(2l_0)^{-1}, 1/2], \ l_0>1$  we apply \refm[III] in condition~$(III)$. By \refm[phieq] and \refm[III] we have for $n$ large enough,
\be |\phi_n(\alpha)|\le \delta_n^{1+\frac{\rho_r}{2}+\epsilon },\ \epsilon>0.\la{boue}\eu
Comparing  the bounds \refm[bou1], \refm[logbou], \refm[boue]
with the asymptotics \refm[I1sim], \refm[basym]
proves \refm[relat].
\hfill\qed

\section{A representation of $W(q)$}\label{Wexp}
In this appendix we derive representations of the function
$W(q)$ in the case of rational  $q>1$.
The infinite product
\be F(z):= \prod_{k=1}^\infty \left(1+\frac{z}{k^q}\right),\quad z\in \BC,\quad q>1,  \la{weir}\eu
is a Weierstrass representation of an entire function $F$ with zeroes
at $\{-k^q, \ k=1,2,\ldots\}$. This follows from Theorem 5.12 in \cite{conw}, since $\sum_{k=1}^\infty k^{-q}<\infty,\ q>1.$ Note that  $W(q)=f(1)=F(1),\ q>1.$
We now show that in the case when $q>1$ is a rational number, a  modification of the argument in \cite{watwhit}, p.238 allows us to
decompose the value $F(1)$ in \refm[weir] into a finite product of values of a canonic entire function of finite rank. (For the definition of a rank of entire function see Chapter $XI$ in \cite{conw}).
Let $q=\frac{m_1}{m_2},  $
where $m_1> m_2\ge 1$ are co-prime integers.
 We write $$1+k^{-\frac{m_1}{m_2}}=\prod_{l=1}^{m_1}\frac{k^{\frac{1}{m_2}}-\alpha_l(m_1)}{k^{\frac{1}{m_2}}}=\prod_{l=1}^{m_1} \left(1-\frac{\alpha_l(m_1)}{k^{\frac{1}{m_2}}}\right),$$
where  $$\alpha_l(m_1)=
\exp\left(\frac{\pi(2l-1)}{m_1}i\right),\quad l=1,\ldots,m_1$$
are all $m_1$-th roots of $-1$, such that $0<\arg(  \alpha_l(m_1))<2\pi, \ l=1,\dots ,m_1.$

Consequently,  \be W(q)=\prod_{k=1}^\infty \prod_{l=1}^{m_1} \left(1-\frac{\alpha_l(m_1)}{k^{\frac{1}{m_2}}}\right),\ \quad q=\frac{m_1}{m_2}.\la{pjs}\eu
%$$ \prod_{l=1}^{m_1} \prod_{k=1}^\infty \big(1-\frac{\alpha_l(m_1)}{k^{\frac{1}{m_2}}}exp(\sum_{l=1}^{m_1-1}\big) .$$
%\be \prod_{l=1}^{m_1} \prod_{k=1}^\infty \Big(\big(1-\frac{\alpha_l(m_1)}{k}\big) %e^{\frac{\alpha_l(m_1)}{k}}\big)\Big).\la{hjuk}\eu
Next, introduce  the function \be \tilde{f}(z):=\prod_{k=1}^\infty\left(1+ \frac{z}{k^{\frac{1}{m_2}}}\right)\exp\left(\sum_{p=1}^{m_2}\frac{(-z)^p}{k^{\frac{p}{m_2}}p}\right),\quad z\in \BC, \la{ggam}\eu
which is a canonical form of an entire function of finite rank $m_2$ with zeroes
$$
\left\{-k^{\frac{1}{m_2}}, \quad k=1,2,\ldots\right\}.
$$
Observing that $\sum_{l=1}^{m_1}\left(\alpha_l(m_1)\right)^p=0, \ p=1,\ldots,m_2,$ by the definition of $\alpha_l(m_1)$, $l=1,\dots, m_1$, and rewriting \refm[ggam] as
$$\tilde{f}(z):=\prod_{k=1}^\infty\left(1+ \frac{z}{k^{\frac{1}{m_2}}}\right)\prod_{p=1}^{m_2}\exp\left(\frac{(-z)^p}{k^{\frac{p}{m_2}}p}\right)$$
we  derive from \refm[pjs]:
\be W(q)=\prod_{l=1}^{m_1} \tilde{f}(-\alpha_l(m_1)),\la{vad}\eu
for rational $q>1.$

For $m_2>1,$ we will consider now  the function
\be \tilde{\Gamma}(z):= e^{Q(z)}\frac{1}{z\tilde{f}{(z)}},\la{chb}\eu
where $Q(z)$ is a polynomial in $z$ that will be defined below.
The preceding discussion yields that $\tilde{\Gamma} $ is a meromorphic function in $\BC$ with   simple poles at
$(-k^{\frac{1}{m_2}}),\ k=0,1,\ldots$.
Now our purpose will be to obtain for the function $\tilde{\Gamma} $ an analog of Gauss formula for the gamma function. We recall the definition of
generalized Euler constants:
\begin{eqnarray*}
\gamma_\alpha&=&\lim_{n\to \infty} \big(\sum_{k=1}^n k^{-\alpha}- \int_1^n x^{-\alpha}dx \big)\\
&=&
\lim_{n\to \infty}
\left\{
                    \begin{array}{ll}
                      \sum_{k=1}^n \frac{1}{k}-\log n, & \hbox{if  }\ \alpha=1; \\
                      \sum_{k=1}^n \frac{1}{k^{\alpha}}-\frac{n^{1-\alpha}-1}{1-\alpha} , & \hbox{if  }\ 0<\alpha<1.
                    \end{array}
                  \right.
\end{eqnarray*}
(Note that $\gamma_1=\gamma$ is the standard Euler constant).
This allows us to write the function $\frac{1}{z\tilde f(z)},\ z\in \BC$ in the following form:
\begin{eqnarray*}
\frac{1}{z\tilde f(z)}&=&
\lim_{n\to \infty} \frac{1}{z}\prod_{k=1}^n\left(\frac{k^{\frac{1}{m_2}}}{z+k^{\frac{1}{m_2}}} \right)\exp\left(-\sum_{p=1}^{m_2}\frac{(-z)^p}{k^{\frac{p}{m_2}}p}\right)\\
&=&
\exp{\Big(-\sum_{p=1}^{m_2}(-1)^p\frac{z^p}{p}\gamma_{(p/m_2)}}\Big)\\
&&\times
\lim_{n\to \infty}\frac{(n!)^{\frac{1}{m_2}}}{\prod_{k=0}^n\left(z+k^{\frac{1}{m_2}}\right)}
n^{\frac{-(-z)^{m_2}}{m_2}}
\exp\left(-m_2\sum_{p=1}^{m_2-1}
\frac{(-1)^pz^p\left(n^{\frac{m_2-p}{m_2}}-1\right)}
{p(m_2-p)}
\right).
\end{eqnarray*}
Setting now in \refm[chb] $Q(z)= \sum_{p=1}^{m_2}(-1)^p\frac{z^p}{p}\gamma_{(p/m_2)},$
we arrive at the desired representation of the function
\be \tilde{\Gamma}(z)= \lim_{n\to \infty}\frac{(n!)^{\frac{1}{m_2}}}{\prod_{k=0}^n\left(z+k^{\frac{1}{m_2}}\right)}
n^{\frac{-(-z)^{m_2}}{m_2}}\exp\left(-m_2\sum_{p=1}^{m_2-1}\frac{(-1)^pz^p\left(n^{\frac{m_2-p}{m_2}}-1\right)}{p(m_2-p)}
\right). \la{dcuk}\eu
Under $m_2=1$, \refm[dcuk] becomes the Gauss formula for the gamma function.

In the case  $q>1$ is an integer,
\refm[vad] conforms to the explicit expression for $W(q)$ in \cite{watwhit}, p.238-239.
In fact,
after substituting in \refm[ggam] $m_2=1$ and $p=1$ we have
$$\tilde{f}(z)= \prod_{k=1}^\infty\left(1+ \frac{z}{k}\right)e^{-\frac{z}{k}},
\quad z\in \BC $$
and by the Weierstrass factorization theorem for the Gamma function,
$$\tilde{f}(z)=\frac{e^{-\gamma z}}{\Gamma(1+z)},\quad z\in \BC\backslash \{-1,-2,...  \}.$$
where $\gamma$ is Euler's constant.
Thus, when $q>1$ is an integer,
$$
W(q)=\prod_{l=1}^q \Big(\Gamma(1-\alpha_l(q))\Big)^{-1}.
$$
Taking into account that the numbers $\alpha_l(q),\ l=1,\ldots,q$ are pairwise conjugate and that $\Gamma(\bar{z})=\overline{\Gamma(z)},\ z\in \BC,$ the last expression can be written as follows:
$$
W(q)=\prod_{l=1}^{[q/2]} \Big(\vert\Gamma(1-\alpha_l(q))\vert^2\Big)^{-1},
\quad q>1.
$$

\end{appendices}

\subsection*{Acknowledgements}
Our paper was benefited by the remarks of an anonymous referee who, among other things, pointed an error in the original proof of the case $q>1$ in Theorem~\ref{weighted}.

\end{document}